\documentclass[draftcls,onecolumn, 11pt]{IEEEtran}

\usepackage[latin1]{inputenc}
\usepackage[T1]{fontenc}
\usepackage{amssymb}
\usepackage{amsmath}
\usepackage{dsfont, graphicx}
\usepackage{cite}

\newtheorem{algo}{Algorithm}

\def\C{\mathbb{C}}

\def\1{\mathds{1}}

\title{OMP-type Algorithm with Structured Sparsity Patterns for Multipath Radar Signals}
\author{Tabea~Rebafka,   C\'eline~L\'evy-Leduc     and~Maurice~Charbit,~\IEEEmembership{Member,~IEEE.}
\thanks{T.~Rebafka is with the LPMA at UPMC, University Paris 6, France.
 e-mail: tabea.rebafka@upmc.fr.}%
 \thanks{C. Lévy-Leduc and M. Charbit  are with the CNRS/LTCI of Telecom ParisTech, France.
 e-mail: \{celine.levy-leduc, maurice.charbit\}@telecom-paristech.fr.}
 \thanks{This work was financially  supported by the French National Agency: Direction Générale de l'Armement (DGA).}
 }

\begin{document}

\maketitle

\begin{abstract}
A transmitted, unknown radar signal is observed at the receiver through more than one path in additive noise. The aim is to recover the waveform of the intercepted signal and to simultaneously estimate the direction of arrival (DOA). We propose an approach exploiting the parsimonious time-frequency representation of the signal by applying a new  OMP-type algorithm for structured sparsity patterns. An important issue is the scalability of the proposed algorithm since high-dimensional models shall be used for radar signals.  Monte-Carlo simulations for modulated signals illustrate the good performance of the method even for low signal-to-noise ratios and a gain of 20 dB for the DOA estimation compared to some elementary method.
\end{abstract}

%

\section{Introduction}
The aim of SIGnals INTelligence (SIGINT) is to intercept as much  information as
possible on received signals. For radar sources, the
parameters of interest are the Direction Of Arrival (DOA) of the
direct path (assuming it is observed), but also, if possible, the
carrier frequency and the modulation scheme. In
practice,  due to the presence of multipath propagation, 
the interception of radar sources remains a difficult problem. 
Indeed, multipath propagation creates nuisance parameters that have to
be taken into account such as DOAs and relative Times of Arrival (TOA) of reflected paths, 
relative signal power levels on the different paths and additive noise variance.

In this paper, we shall consider a unique narrowband farfield source
propagating through several paths, corrupted with  additive Gaussian
white noise and observed on a multi-sensor array.
This issue is actually very involved since the observed signals
on the different sensors are highly correlated and the different
replicates are likely to overlap. This prevents us from using methods such as
MUSIC, ESPRIT or MVDR which do not work very well with coherent
signals \cite{viberg1996}. These approaches, based on subspace methods
or on maximum likelihood estimation, have been extended  by \cite{zhang2009}, \cite{vanderveen1997},
\cite{weng1996} and \cite{shan1985} to deal with
coherent sources. However, to the best of our knowledge, estimating
the waveform, which is an infinite-dimensional parameter, in a
nonparametric way has never been considered for multipath signals. To solve this problem, our
main idea consists in exploiting the sparsity of the time-frequency
representation of the waveform, 
which means that the selection of possibly useful signal components from a huge collection of candidate signals
may be driven by a sparsity constraint. 
Such kind of approach 
was initially proposed by \cite{fuchs99}  but only for the estimation of the time-delays, which is a finite dimensional parameter, and 
is extended here for estimating the DOAs and TOAs  as well as the waveform.

For dealing with the estimation in sparse linear regression models, the
Lasso \cite{tibshirani96} and the greedy orthogonal matching pursuit (OMP)   
(\cite{pati93}, \cite{mallat94}) have become very popular tools.
In our case, an inspection of our model shows that the parameter 
vector is not sparse in an arbitrary way:
its sparsity pattern has a specific structure. In order to include
prior information concerning the sparsity structure, different approaches
have  recently been proposed. On the one hand, there are methods based on
composed $\ell_1/\ell_2$-penalties  (elastic nets \cite{zou05}; 
fused Lasso \cite{tibshirani05}; group Lasso \cite{yuan06}; 
composite absolute penalty \cite{zhao09}; overlapping groups
\cite{jacob09,jenatton09}), on the other hand, \cite{lozano09}
proposed the group-OMP method for non-overlapping groups, which is an extension of the orthogonal matching pursuit 
to structured solutions.  

In this paper we explain how to extend the
group-OMP to overlapping groups and how to deal with the
scalability of our algorithm which is a crucial issue in the context of
intercepted  radar signals since the model dimension is usually very high.
A simulation study shows that  the proposed method performs well for
composed and modulated signals even when the signal-to-noise ratio is
low. A gain of more than 20 dB for the DOA estimation is achieved in comparison with some elementary method.


The paper is organized as follows. In Section \ref{sec_model} the mathematical model for radar signals is introduced, together with a reformulation of the problem in the form of a sparse, partly linear model. In Section \ref{sec_algo} an OMP-type algorithm is developed which is well adapted to our model. Section \ref{sec_scalable} deals with the scalability issue of the algorithm, while Section \ref{sec_stop} presents an appropriate stopping criterion for the OMP-type algorithm. Finally, the  experimental results of Section \ref{sec_simul} illustrate the performance of the new algorithm and a comparison with some elementary method is provided. Conclusions are made in the final  Section \ref{sec_concl}.

\section{Model}\label{sec_model}

\subsection{Signal Modelling}
Let us denote by $s_0(t)$ the original radar signal emitted by the source, where $t$ is the time. Due to  propagation  and  reflexion on obstacles, the  signal arrives via several pathways on the receiver. With each  pathway we associate a direction of arrival, say $D_u$, a time-delay, say $T_u$ and a complex constant $A_u$ representing the attenuation  of the  signal on the $u$-th pathway. The signal is detected by an antenna with $C$ sensors whose array response is denoted by $\vec a(\cdot)$. Then the signal arriving at the sensor array at time $t$ is a $C$-dimensional vector given by
\begin{equation}\label{def_signal_antenna}
	\vec s(t)=\sum_{u=1}^U A_u\,s_0(t-T_u)\, \vec a(D_u)\;,
\end{equation}
where $U$ denotes the number of propagation pathways, which is unknown. As the time lags $T_u$ are typically much shorter than the total length  of the signal $s_0$, the arriving signal $\vec s$ is the superposition of several delayed replicates  $s_0$ with different amplitudes.

In practice, the signal $\vec s$ is observed with some additional noise, say $\vec\varepsilon$, and at a sample frequency, say $F_s$. With $\Delta_s = 1/F_s$ the observations are hence given by
\begin{equation}\label{def_observed_signal}
 	\vec y_m = \vec s(m\Delta_s) + \vec \varepsilon_m\;,\qquad m=1,\dots,M\;,
\end{equation}
where $\vec \varepsilon_m\in\C^C$ is a noise vector. Generally one considers   circular gaussian independent random vectors with zero mean and covariance $\sigma^2I_C$.

\subsection{Linear Regression Model}
In practice, all  model  parameters   are unknown, except the array response function $\vec a(\cdot)$, which is determined by the array geometry or by calibration. As the relationship between the parameters is rather  involved, we propose to reformulate the estimation problem by  linearizing the model. This is obtained, on the one hand, by representing the waveform in an overcomplete basis and, on the other hand, by discretizing  the parameter spaces of the other parameters. This leads to a linear model with a high-dimensional parameter vector and nonlinear constraints. As this vector is structured and sparse, \textit{i.e.} most entries are zero,  regularization methods with structured sparsity inducing norms can be applied \cite{zhao09} .

\subsubsection{Decomposition of the signal }\label{subsubsec_decomposition}
To estimate the waveform, the radar signal $s_0$  is represented in some given (overcomplete) basis or dictionary $\mathcal{D} =\{\varphi_j, j=1,\dots,J\}$, such that
\begin{equation}\label{def_signal_basis}
	s_0 = \sum_{j=1}^J\beta_j\varphi_j\;,
\end{equation}
with appropriate coefficients $\beta_j\in\C$. Hence, estimating   $s_0$ becomes recovering the coefficients $\beta_j$. When using appropriate, large dictionaries, a small number of elements $\varphi_j$ is sufficient to reconstruct $s_0$. In other words, it is assumed  that the dictionary is such that most coefficients $\beta_j$ are   zero and hence that the representation of $s_0$ is {sparse}.

\subsubsection{Propagation model}
To estimate the parameters $A_u$, $D_u$ and  $T_u$ by discretization, we introduce  grids  $\{\tau_1,\dots,\tau_P\}$ and $\{\theta_1,\dots,\theta_Q\}$ of potential values of the delay times $T_u$ and the angles of arrival $D_u$ respectively, as in \cite{fuchs99,malioutov05}. Then $\vec s(t)$ can be rewritten as
\begin{equation}\label{def_antenna_signal_discrete}
	\vec s(t)=\sum_{p=1}^P\sum_{q=1}^Q \alpha_{p,q} ~s_0(t-\tau_p)\, \vec a(\theta_q)\;,
\end{equation}
$$\text{with }\quad\alpha_{p,q} = \sum_{u=1}^UA_u\1\{T_u=\tau_p, D_u=\theta_q\}\;,$$
where $\1$ is the indicator function. Clearly, there are  only $U$ coefficients $\alpha_{p,q}$ that are nonzero. Hence, when fine grids are used, \textit{i.e.}  $P$ and $Q$ are large,  this is another sparse estimation problem. As the number of pathways $U$ is generally unknown, model selection is required, \textit{i.e.} the estimation of the model order $U$.

By combining (\ref{def_signal_basis}) and  (\ref{def_antenna_signal_discrete}), the observation $\vec y_m$ reads
\begin{equation}\label{def_observed_signal_discrete}
	\vec y_m =  \sum_{j=1}^J\sum_{p=1}^P\sum_{q=1}^Q\alpha_{p,q} \beta_j\varphi_j(m\Delta_s-\tau_p)\, \vec a(\theta_q) + \vec \varepsilon_m\;,
\end{equation}
for $m=1,\dots,M$. It is convenient to represent the model in matrix notation. Let  $\mathbf{X}$ be the  matrix of size $CM\times JPQ$, where all known quantities are stored. More precisely, the $(j,p,q)$-th predictor, \textit{i.e.} the column ${x}_{j,p,q}$ of $\mathbf{X}$, is given by
$$x_{j,p,q} =\left(\begin{array}{c} \varphi_j(\Delta_s-\tau_p)\,\vec a(\theta_q)\\
    \vdots\\
    \varphi_j(M\Delta_s-\tau_p)\,\vec a(\theta_q)
    \end{array}\right)  \in\C^{CM} \;.
$$
In other words, ${x}_{j,p,q}$ corresponds to the signal received by the  antenna (without noise) when the emitted signal is the element $\varphi_j$ arriving with delay $\tau_p$ and angle $\theta_q$. Denote the observation vector by $Y = (y_1^T,\dots,y_m^T)^T\in\C^{CM}$   and the noise vector by $\varepsilon = (\vec \varepsilon_1^T,\dots,\vec \varepsilon_m^T)^T\in\C^{CM}$. Finally, let the Kronecker product $\beta \otimes\alpha$ of the vectors $\beta = (\beta_1,\dots,\beta_J)^T$ and $\alpha = (\alpha_{1,1},\dots, \alpha_{P,Q})^T$ be the $JPQ$-vector $w$ defined as
\begin{equation}\label{rel_coeff_w_al_beta}
w=\left(\begin{array}{c}\vdots\\
w_{j,p,q}\\\vdots\end{array}\right)
= \left(\begin{array}{c}
\beta_1\alpha_{1,1}\\
\beta_1\alpha_{1,2}\\
\vdots\\
\beta_j\alpha_{p,q}\\
\vdots\\
\beta_J\alpha_{P,Q}
\end{array}\right)\:.
\end{equation}
Then model (\ref{def_observed_signal_discrete}) can be written as
\begin{equation}\label{def_model_al_be}
	Y= \mathbf{X}(\beta\otimes\alpha)+ {\varepsilon}\quad\text{ with }  \alpha\in\C^{PQ}, \beta\in\C^J\;.
\end{equation}
Finally, introduce the following set
\begin{align*}
	\mathcal{W}  &= \{w\in\C^{JPQ}: w = \beta \otimes\alpha,\; \alpha \in \C^{PQ}, \beta\in\C^J \}\;, \label{def_W_searchspace}
\end{align*}
which is a subspace of $\C^{JPQ}$  that can be parametrized using $J+PQ$ variables. Remark that, given a vector $w\in\mathcal{W}$, one can recover the	vectors $\alpha$ and $\beta$ only up to a multiplicative constant,	which does not matter in our application. Thus, model
(\ref{def_observed_signal_discrete}) can be written as the following linear regression model
\begin{equation}\label{def_lin_model}
	Y= \mathbf{X}{w}+ {\varepsilon}\quad\text{ with }  w\in\mathcal{W},
\end{equation}
where the design matrix $\mathbf{X}$ is known and the parameter vector $w$ is to be estimated from the observations $Y$.  The particularity of this regression model is that the search space for $w$ is not the entire space $\C^{JPQ}$ but the smaller, nonconvex
set $\mathcal{W}$. 

\section{Estimation procedure}\label{sec_algo}
The problem now is  to estimate the vectors $\alpha$ and $\beta$ by fitting model (\ref{def_model_al_be}), or  equivalently, to estimate $w$ in model (\ref{def_lin_model}). The advantage of (\ref{def_lin_model}) over (\ref{def_model_al_be}) is the linearity of the model.
Since the dimension of the space $\mathcal{W}$  is much higher than that  of the observations $Y$,  the minimization problem
\begin{equation}\label{def_min_pb_w}
	\min_{w\in\mathcal W} \|Y-\mathbf X w\|^2
\end{equation}
is  ill-posed. However, we note that $w$ must be a sparse vector, since $\alpha$ and $\beta$ are assumed sparse. It is well known that sparsity leads to stable estimation procedures based on a regularized fitting criterion with  an $\ell_1$-penalty (Lasso \cite{tibshirani96}; basis pursuit \cite{chen98}; Lars \cite{efron04}) or greedy algorithms (as matching pursuit and orthogonal matching pursuit (OMP) \cite{pati93,chen94}). A closer inspection of $\mathcal W$ reveals some structure in the sparsity patterns, that means, there are constraints on the distribution of the zero entries for all vectors $w$ in $\mathcal W$. More precisely, the set of indices $\mathcal S=\{(j,p,q): \text{ such that } w_{j,p,q}=0\}$ is called 
the sparsity pattern of some given vector $w$. Indeed,  all elements $w$ of $\mathcal W$ verify that $\alpha_{p',q'}=0$ implies that the components $w_{j,p',q'}=0$ for all $j$, and likewise, if $\beta_{j'}=0$ then $w_{j',p,q}=0$ for all $p$ and $q$. Thus, $w$ must not be
sparse in an arbitrary way, but its sparsity pattern has to respect some structure.  Modern estimation procedures  force specific structures of the sparse vector. On the one hand, there are methods based on  composed $\ell_1/\ell_2$-penalties  (elastic nets \cite{zou05}; fused Lasso \cite{tibshirani05}; group Lasso \cite{yuan06}; composite absolute penalty \cite{zhao09}; overlapping groups \cite{jacob09,jenatton09}), on the other hand,  the orthogonal matching pursuit can be extended to structured solutions   (group-OMP for non-overlapping groups \cite{lozano09}). 

Our problem can be related to a penalized minimization problem of the form
\begin{align}\label{def_pen_pb}
	\min_{w\in\C^{JPQ}}\left\{ \|Y-\mathbf X w\|^2+\Omega(w)\right\}\;,
\end{align}
with a structured $\ell_1/\ell_2$-penalty $\Omega(w)$ (detailed below) in the sense that the  sparsity pattern of a solution of this problem, say $\tilde w$, respects the required structure. However, there is no guarantee that the nonzero components of  $\tilde w$ satisfy relation (\ref{rel_coeff_w_al_beta}), that is, that there are  coefficients $\alpha_{p,q}$ and $\beta_j$ such that $\tilde w_{j,p,q}= \alpha_{p,q}\beta_j$  and hence $\tilde w$ may not belong to $\mathcal W$. The problem lies in the nonzero components of $\tilde w$.  To compute the solution of (\ref{def_pen_pb}),  different solutions have been proposed in the literature as the active set algorithm in in \cite{jenatton09b} and \cite{jenatton09}. In this paper we propose  an extension of the group-OMP \cite{lozano09} for overlapping groups.  Now, to fit the model given in (\ref{def_model_al_be})  we propose the following procedure consisting of two steps, that are developed below in more detail.

\noindent \textit{Step 1.  (Model selection)} 
Solve problem (\ref{def_pen_pb}) via an OMP-type algorithm to obtain a solution $\tilde w$ with admissible sparsity structure. Identify the nonzero components $\tilde w_{j,p,q}$ and the indices $\mathcal A$ of the associated coefficients $\tilde\alpha_{p,q}$ and $\tilde\beta_j$ that must be nonzero.

\noindent \textit{Step 2. (Estimation in reduced dimension)}
Reduce the dimension of the regression matrix $\bf X$ by keeping the  predictor $x_{j,p,q}$ only if $\tilde w_{j,p,q}\neq0$.   Then compute the  least squares estimator in model (\ref{def_model_al_be}) by the Nelder-Mead simplex method. This problem is now well-posed because of the largely reduced dimension.

\subsection{Model selection step} 
To specify the penalty $\Omega(w)$ in (\ref{def_pen_pb}), we  describe the structure of the sparsity pattern $\mathcal S$ of $w\in\mathcal W$. To this end, for a given couple $(p,q)$, we denote the set of indices $G_{p,q}^\alpha  =\left\{  (j,p,q), ~j\in\{1,\dots,J\} \right\}$,  and likewise, for a given $j$, the set of indices $G_{j}^\beta = \left\{ (j,p,q), ~p\in\{1,\dots,P\}, q\in\{1,\dots,Q\}\right\}$.  It is clear that the sparsity pattern $\mathcal S$ of any $w$ in $\mathcal W$ can be written as the union of all sets  $G_{p,q}^\alpha$ and $G_{j}^\beta$ with $\alpha_{p,q}=0$ and $\beta_j=0$.  

Given a set of indices $G$, let $\|w_G\|_2$ denote the $\ell_2$-norm of the subvector of $w$ composed of the elements $w_{j,p,q}$ with $(j,p,q)\in G$. Then one can show that the solution, say $\tilde w$, of (\ref{def_pen_pb})  with the following penalty
\begin{equation}\label{eq_penalized_pb}
	\Omega(w) =\lambda_{1}\sum_{p,q}\|w_{G_{p,q}^\alpha}\|_{2} +\lambda_{2}\sum_{j}\|w_{G_{j}^\beta}\|_{2} \;,
\end{equation}
where $\lambda_{1}, \lambda_{2}>0$ are regularization parameters,  has the required sparsity structure, \textit{i.e.} the sparsity pattern of $\tilde w$ can be written as the union of sets $G_{p,q}^\alpha$ and $G_{j}^\beta$. Note that these groups are  possibly overlapping. Indeed, the penalty can be viewed as  a mixture of $\ell_{1}$- and $\ell_{2}$-norms appropriate for overlapping group structures. More precisely, the penalty is the sum (\textit{i.e.} the $\ell_{1}$-norm) of the $\ell_{2}$-norm of subvectors $w_{G}$. Likewise to the standard Lasso, the $\ell_{1}$-norm entails that for several groups $G$ the  terms $\|w_{G}\|_{2}$ are zero, implying in turn that all entries of $w_{G}$ are exactly zero. By using two  regularization parameters $\lambda_{1}>0$ and $\lambda_{2}>0$  one can use  different degrees of sparsity for the vectors $\alpha$ and $\beta$. 

To solve the problem (\ref{def_pen_pb}) with penalty $\Omega(w)$ given by (\ref{eq_penalized_pb}) we propose an OMP-type algorithm. To describe the algorithm we introduce  the collection  $\mathcal G=\{G_g, g\in \mathcal I\}$ of all sets  $G_{p,q}^\alpha$ and $G_{j}^\beta$ indexed by a common index $g$ with index set $\mathcal I$.  Note that the cardinality of $\mathcal I$ is $PQ+J$. For a set of indices $A$, denote by $\mathbf X_A$ the matrix made of the corresponding predictors $x_{j,p,q}$ of $\mathbf X$ such that $(j,p,q)\in A$. 

Recall that the principle of OMP consists in  adding iteratively the predictor to the current solution $w^{(t)}$ which is the most correlated with the current residual $r^{(t)}=Y-\mathbf Xw^{(t)}$. In \cite{lozano09} this principle is extended to selecting non-overlapping groups of variables, where $\mathcal G$ is a set of pairwise disjoint sets $G_g$ and  one adds all variables of group $G_{g^*}$ if $\|\mathbf X_{G_{g^*}}^Tr^{(t)}\|_2^2=\max_{g\in\mathcal I}\|\mathbf X_{G_g}^Tr^{(t)}\|_2^2$.  Proceeding in this way guarantees that at every step of the algorithm, the current solution respects the required sparsity structure, which means that the sparsity pattern of the current solution $w^{(t)}$ equals the union of some groups $G_g\in\mathcal G$.  

Now this procedure can be extended to our context with overlapping group structures, where the aim remains the same: the sparsity pattern of every intermediate solution  $w^{(t)}$ must be the union of  some groups $G_g\in\mathcal G$.  The difference to the group-OMP of \cite{lozano09} is at the level of the selection of variables that may enter the solution. Indeed, due to the overlapping group structure, the sets of  variables that may be activated at the next step are not simply the groups $G_g\in\mathcal G$.  To be more precise on this issue, we introduce the set $\mathcal Z^{(t)}\subset\mathcal I$ of indices of zero groups $G_g$ associated with the current solution $w^{(t)}$, such that  $\mathcal S^{(t)}=\cup_{g\in\mathcal Z^{(t)}} G_g$ is the sparsity pattern of $w^{(t)}$.  Clearly, the indices of the non zero entries of $w^{(t)}$ are given by $\mathcal A^{(t)}=(\mathcal S^{(t)})^c$, where $A^{c}$ denotes the complementary set of $A$ with respect to $\mathcal G$. Now the candidate groups of components that may be activated at the next step are obtained by deleting one of the indices from the zero groups $\mathcal Z^{(t)}$. That means, every potential group of indices has the form $\mathcal P_{g_0} = (\cup_{g\in\mathcal Z^{(t)}\backslash \{g_0\}} G_g)^c\backslash \mathcal A^{(t)}$ where $g_0\in\mathcal Z^{(t)}$. Then, following the philosophy of OMP, one activates the components given by $\mathcal P_{g^*}$ if $\|\mathbf X_{\mathcal P_{g^*}}^Tr^{(t)}\|_2^2=\max_{g_0\in\mathcal Z^{(t)}}\|\mathbf X_{\mathcal P_{g_0}}^Tr^{(t)}\|_2^2$.   Once a group $\mathcal P_{g^*}$ is selected, one derives the new set of active indices $\mathcal A^{(t+1)} = \mathcal A^{(t)}\cup \mathcal P_{g^*}$ and updates  the current solution by solving the least squares problem of reduced dimension $w_{\mathcal A^{(t+1)}}^{(t+1)}=\arg\min_{w} \| Y-\mathbf X_{\mathcal A^{(t+1)}}w\|_2^2$.

We notice that the algorithm can be viewed as a walk through a directed acyclic graph, which also appears in the algorithm proposed in \cite{jenatton09}. Formally, our OMP-type algorithm can be described as follows.

\bigskip

\begin{algo}[OMP-type algorithm for overlapping groups] 
\it Initialization.~~
\rm Set $t=0$.\\
\texttt{ Repeat until} \rm some stopping criterion is satisfied.\\
\indent \texttt{ If} $t=0$ (initial step)\\
\indent \indent \rm Compute $\displaystyle (j^*,p^*,q^*) =\arg\max_{(j,p,q)}\|x_{(j,p,q)}^TY\|_2^2$.\\
\indent \indent Denote by $g^*_\alpha$ and  $g^*_\beta$  the group indices associated with\\
\indent \indent  the groups $G_{p^*,q^*}^\alpha$ and $G_{j^*}^\beta$, respectively\\
\indent \indent Put the set of zero  group indices $\mathcal Z^{(1)}= \mathcal I\backslash \{g^*_\alpha,g^*_\beta\}$.\\
\indent \indent Put the set of active indices $\mathcal A^{(1)}=\{(j^*,p^*,q^*)\}$.
\indent \texttt{  Else}\\
\rm \indent\indent  For each $g_0\in\mathcal Z^{(t)}$ set
$$\mathcal P_{g_0}=\left(\cup_{g\in\mathcal Z^{(t)}\backslash\{g_0\}}G_g\right)^c\backslash \mathcal A^{(t)}\;.$$
\indent \indent Compute $\displaystyle g^* =\arg\max_{g_0\in\mathcal Z^{(t)}}\|(\mathbf X_{\mathcal P_{g_0}})^Tr^{(t)}\|^2_2$.\\
\indent \indent Update  $\mathcal Z^{(t+1)}=\mathcal Z^{(t)}\backslash\{g^*\}$ and $\mathcal A^{(t+1)}=\mathcal A^{(t)}\cup\mathcal P_{g^*}$.\\
\indent \texttt{  End}\\\rm
\indent Update the current solution 
\begin{equation}\label{def_iter_ols}
\tilde w^{(t+1)}_{\mathcal A^{(t+1)}}=\arg\min_{w\in \C^{|\mathcal A^{(t+1)}|}} \| Y - \mathbf X_{\mathcal A^{(t+1)}}w\|^2_2\;.
\end{equation}
\indent Update the  residual $r^{(t+1)}= Y- \mathbf X_{\mathcal A^{(t+1)}}\tilde w^{(t+1)}$.\\
\indent Set $t=t+1$.\\
\texttt{ End} 
\end{algo}
\bigskip

Note that in any case the algorithm stops after a finite number of iterations, more exactly after $J + PQ - 1$ steps. Then all groups are activated and the solution $\tilde w$ is the ordinary least squares estimator.  The important issue of an appropriate stopping criterion for the algorithm is addressed in Section \ref{sec_stop}. Clearly, the aim is to stop the algorithm earlier  to provide a relevant, sparse solution.

We remark that the proposed algorithm is a natural extension  of the orthogonal matching pursuit. Moreover it is similar to the active set algorithm described in \cite{jenatton09}, which provides an interesting mathematical analysis of the optimization problem given by~(\ref{def_pen_pb}). Nevertheless, our OMP-type algorithm provides an approach which is feasible in practice and computationally faster than the active set algorithm. Indeed, in (\ref{def_iter_ols}) the solution is updated by the ordinary least squares estimator that is known explicitly,  whereas \cite{jenatton09}  proposes to solve the \textit{penalized} problem of reduced dimension by a more expensive second-order cone programming. 

\subsection{Solution of the reduced problem}
The solution of (\ref{def_pen_pb}) is a sparse vector, say  $\tilde w$, whose sparsity pattern serves to identify the indices of the  nonzero coefficients $\alpha_{p,q}$ and $\beta_{j}$.  More precisely, if $\tilde w_{j,p,q} \neq 0$, this implies that the associated coefficients $\alpha_{p,q}$ and $\beta_{j}$  are both non zero. Conversely, $\tilde w_{j,p,q} =0$ implies that at least one of the coefficients is 0. 

To obtain estimates of these nonzero coefficients,  we reduce the dimension of the regression matrix $\bf X$ by keeping the  predictor $x_{j,p,q}$ only if $\tilde w_{j,p,q}\neq0$. Denote the resulting matrix by  $\mathbf{X}_\text{red}$. Then compute the least squares estimator by minimizing
\begin{equation}\label{def_min_red_pb}
	L({\tilde\alpha, \tilde\beta}) = \|Y-\mathbf X_\text{red}(\tilde\beta\otimes\tilde\alpha)\|^2\;,
\end{equation}
where $\tilde\alpha$ and $\tilde\beta$ are  the associated $\alpha$- and $\beta$-vectors of reduced dimension $\tilde P\tilde Q$ and $\tilde J$, respectively  This minimization is a well-posed problem as the dimension of $\mathbf{X}_\text{red}$ is much smaller than that of $\mathbf{X}$.

To compute the solution of (\ref{def_min_red_pb}) one can use the fact that the model is linear in $\tilde\alpha$ for fixed $\tilde\beta$ and conversely. That means for fixed $\tilde\beta$ we define the matrix $\mathbf{X(\tilde \beta)}$ with $\tilde P\tilde Q$ columns by
$$\mathbf{X}_\text{red}(\tilde \beta) = \sum_{j=1}^{\tilde J}\tilde \beta_{j} (x_{j, 1,1}^\text{red},\dots,x^\text{red}_{j, {\tilde P},{\tilde Q}})\;.$$
Then the minimum of $\tilde\alpha\mapsto L({\tilde\alpha, \tilde\beta}) $ is the ordinary least squares estimator given by 
$$\tilde\alpha^\text{OLS}(\tilde\beta) = (\mathbf{X}_\text{red}(\tilde \beta) ^{T}\mathbf{X}_\text{red}(\tilde \beta) )^{-1}\mathbf{X}_\text{red}(\tilde \beta) ^{T}Y\;.$$
Furthermore, as $\tilde\alpha$ and $\tilde\beta$ are identifiable only up to some multiplicative constant, one can set $\tilde\beta_{1}=1$, for instance.  Finally, it remains to minimize
$$ (\tilde \beta_{2},\dots,\tilde \beta_{\tilde J})\mapsto L\left(\tilde\alpha^\text{OLS}((1,\tilde \beta_{2},\dots,\tilde \beta_{\tilde J})^{T}),(1,\tilde \beta_{2},\dots,\tilde \beta_{\tilde J})^{T}
\right)$$
which can be done by the Nelder-Mead simplex method. This is feasible since the number of unknown variables, \textit{i.e.} $\tilde J-1$, is quite low.

\section{Scalability}\label{sec_scalable}
To obtain a good representation of the waveform,  a large dictionary shall be used. Likewise, to avoid biased estimators of angles and time delays, we may use fine grids. However, large dictionary sizes $J$ and large grid sizes $PQ$ entail a huge regression matrix $\bf X$ having $JPQ$ columns (see Section \ref{sec_simul} where $JPQ \approx 6.2\,10^9$). This raises computational difficulties, namely concerning the storage of the regression matrix $\bf X$.

A closer look at the OMP-type algorithm  shows that there are no computations involving  the entire regression matrix. Essentially,  the algorithm  considers one by one all possible groups of variables $G\in\mathcal G$ that may enter the solution in the present iteration. For each   group of variables $G$  the correlation of the associated  predictors with the current residual $r^{(t)}=Y-\mathbf X_G \tilde w^{(t)}$ is computed. Hence, instead of the whole matrix $\mathbf X$, only  the predictors of   group $G$ are  required. 

In short, the storage of  $\bf X$ can be avoided by recomputing the required predictors at each iteration. With such a programming there are almost no limits on the size of the regression matrix, and thus almost arbitrary  dictionaries and grids may be used.

To be more precise, let us have a look at the dictionary used in the following simulation study. The dictionary is composed of sinusoids with different frequencies $f$, different lengths $l$ and different starting points $s$. The index $j$ is hence a triple index $(f,l,s)$.  Now instead of storing the whole vector $x_{j,p,q}$, we just use the triple index $(f,l,s)$ to reconstruct the dictionary element $\varphi_{j}$, whenever it is needed. Then further translation and multiplication with the appropriate steering vector yield the required predictor $x_{j,p,q}$. Finally, we only have to store the steering matrix for all possible angles of arrival, which is a matrix of very moderate size $C\times P$. With this knowledge, the predictor $x_{j,p,q}$ is easily reconstructed at any iteration of the algorithm.

In regression problems it is common to work with a normalized regression matrix, that is all columns have mean zero and standard deviation one. To avoid the computation of the normalization constants for every predictor, they may be computed just once and stored for later usage. Note that the normalization constants do not depend neither on the starting point of the dictionary element nor on the delay on the pathway but only on the frequency $f$ and the length $l$ of the dictionary element and the angle of arrival. Hence,  the number of constants to store remains reasonable.

The most time-consuming step in the algorithm is the scan of all predictors (or more precisely,  of all relevant groups of predictors) in every iteration, to identify the one that is the most correlated with the current residual. It is noteworthy that this step can be speeded up by parallelization of the algorithm. Indeed, there is no specific order to visit the predictors and the problem can be split in several independent subtasks.

\section{Stopping criterion}\label{sec_stop}

The OMP-type algorithm requires an appropriate stopping criterion. 
A common approach to conceive a good stopping criterion  is based on cross-validation \cite{zhao09}. Note that in our context,  we do not have to deal with a single, but with two sparse vectors, namely $\alpha$ and $\beta$ and their degrees of sparsity may not be equal.  Hence, cross-validation would be computationally demanding.

Here we propose a stopping criterion that on the one hand takes into account available prior information on the number of propagation pathways (and hence on the number of $\alpha$-variables) and on the other hand has a data-driven component to determine the total number of non zero $\alpha_{p,q}$- and $\beta_j$-coefficients. 

\subsection{Sparsity of $\alpha_{p,q}$-coefficients}
In the context of intercepted radar signals, one generally expects a small number of relevant propagation pathways. Indeed, we are mainly interested in the recovery of the direct path to determine the principal  DOA.  
 It is hence convenient to introduce an upper limit for the number of paths $U_\text{max}$  and to modify the OMP-type algorithm such that at every iteration  the number of activated $\alpha_{p,q}$-coefficients is checked. If their number has achieved $U_\text{max}$, then we stop activating further $\alpha_{p,q}$-coefficients and concentrate on activating only $\beta_{j}$-coefficients. 

In simulations we observed that this procedure yields satisfactory results. The only important assumption for a good performance is that the direct path is not too weak compared to the indirect pathways, since the algorithm activates the components by their importance in the signal.

Note that, if $U_\text{max}$ is smaller than the `true' number of propagation pathways, then the observed signal $Y$ cannot be completely explained. In other words, the solution obtained by the OMP-type algorithm is not the optimal solution of the minimization problem (\ref{def_pen_pb}).


\begin{figure}
	\begin{center}
         \begin{tabular}{cc}
	\includegraphics[scale=.22]{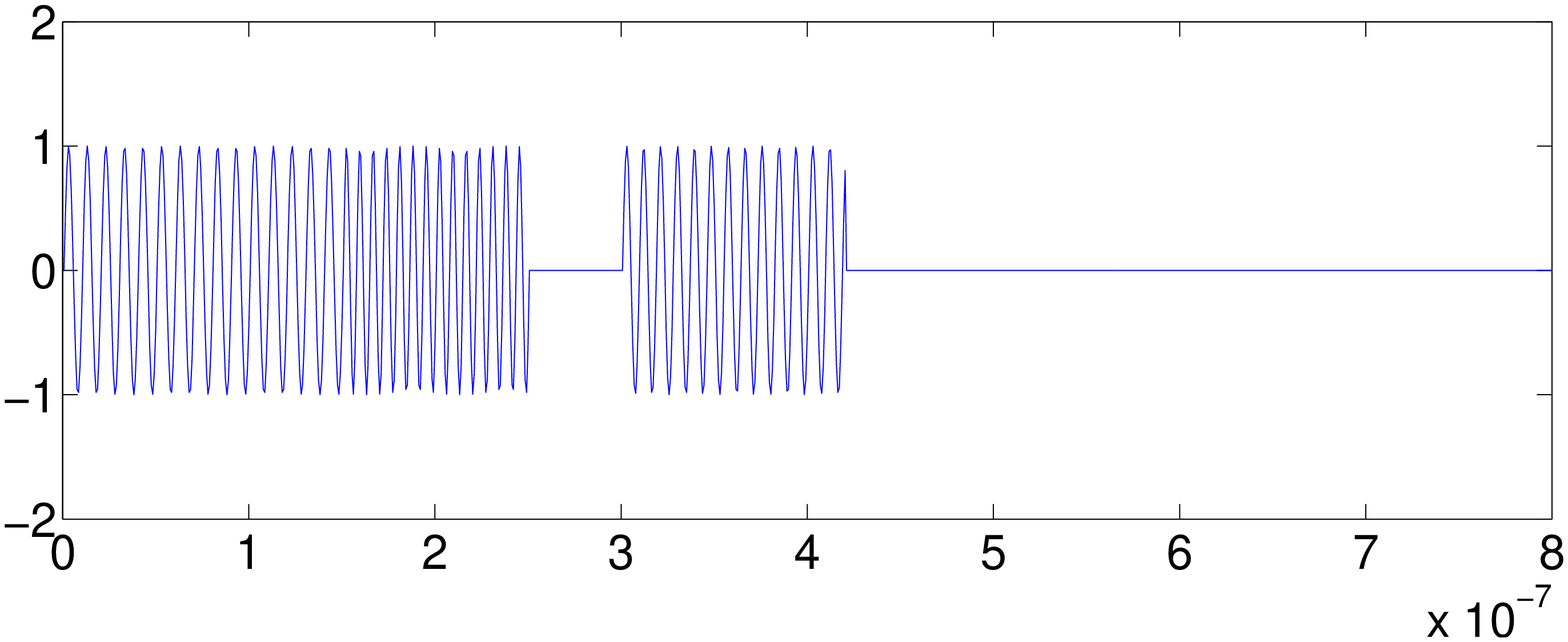}&\includegraphics[scale=.22]{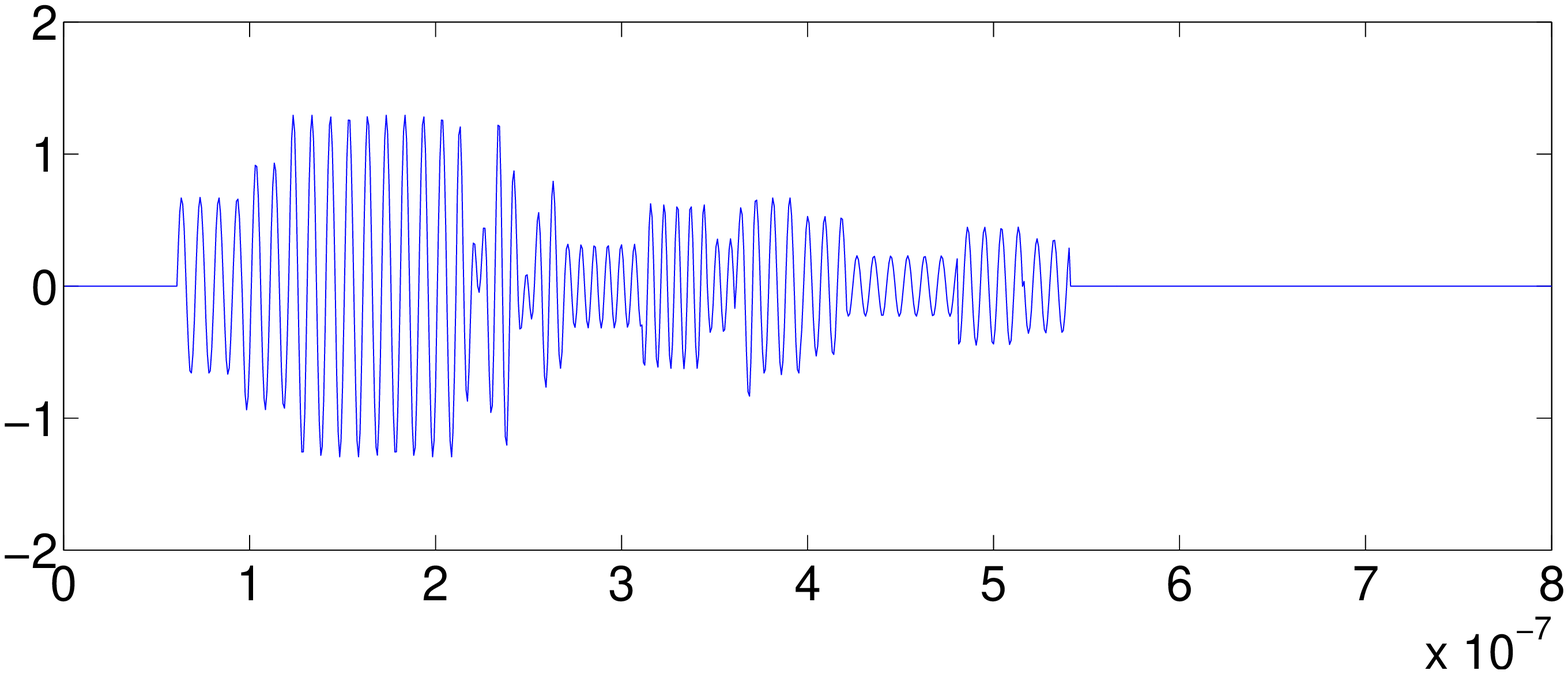}\\
         {\footnotesize	(a) original radar signal}&{\footnotesize (b) multipath signal, no noise, imaginary part}\\
	\includegraphics[scale=.22]{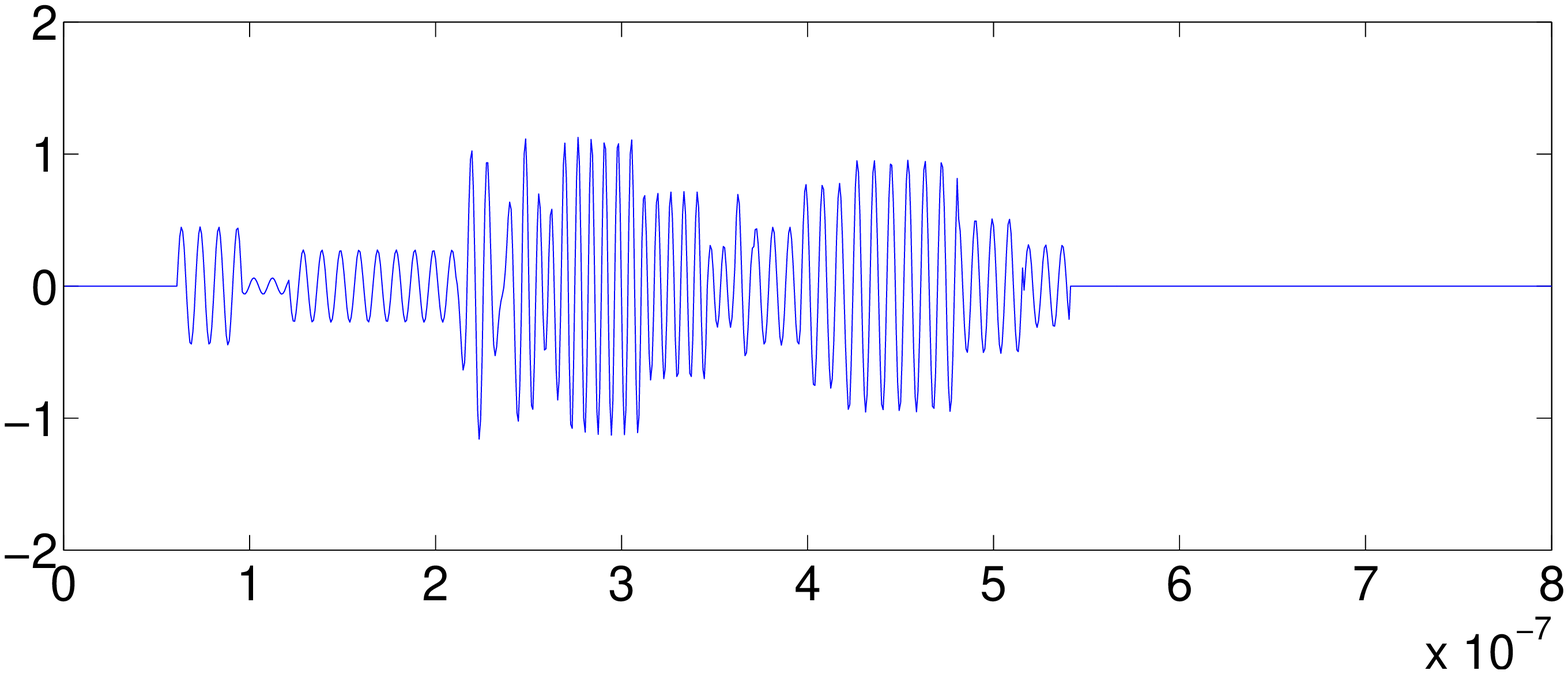}&\includegraphics[scale=.22]{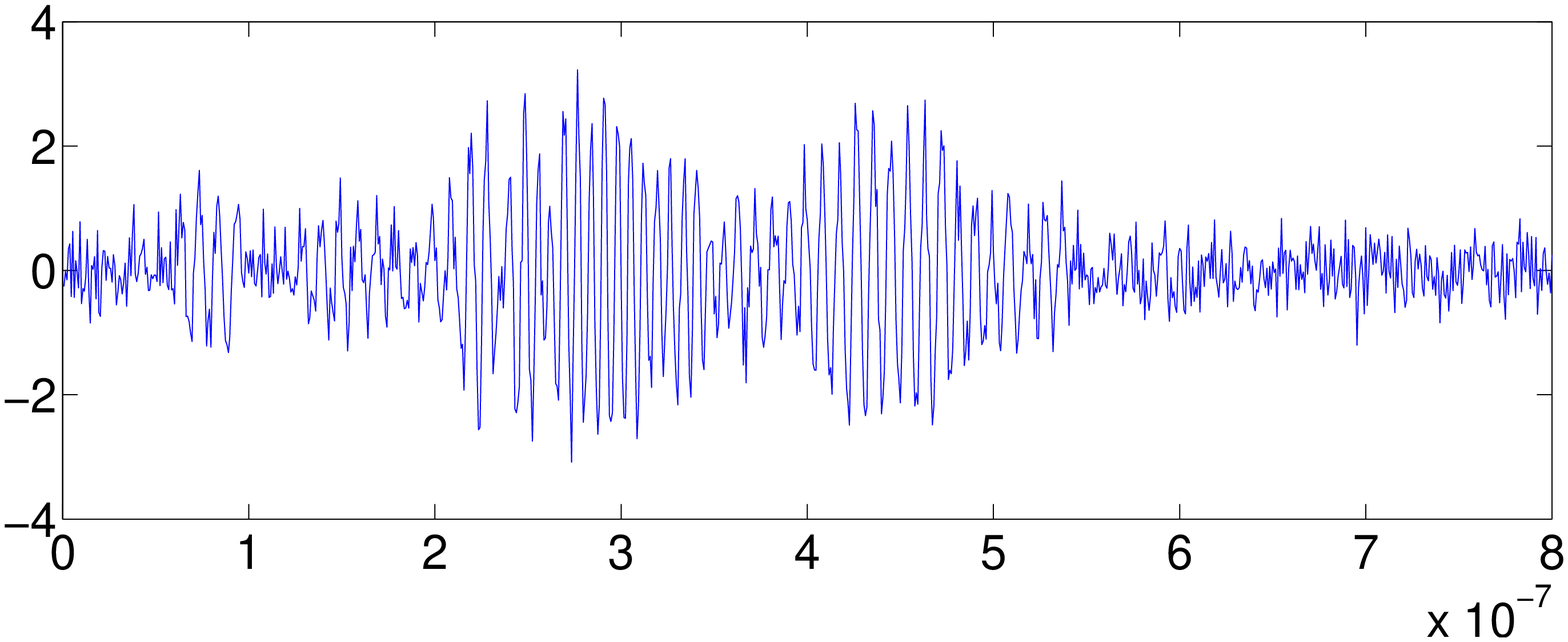}\\
	{\footnotesize (c) multipath signal, no noise,  real part}&{\footnotesize(d) noisy signal, real part}\\
	\includegraphics[scale=.22]{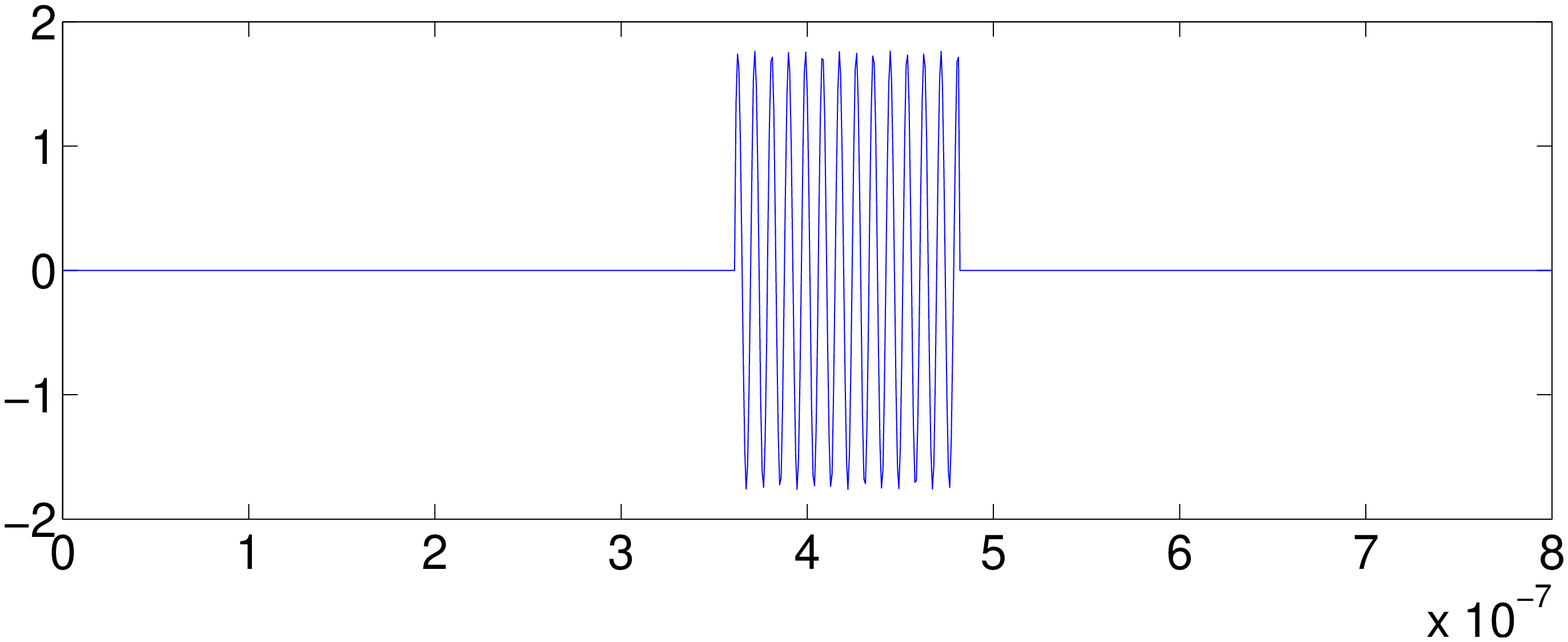}&\includegraphics[scale=.22]{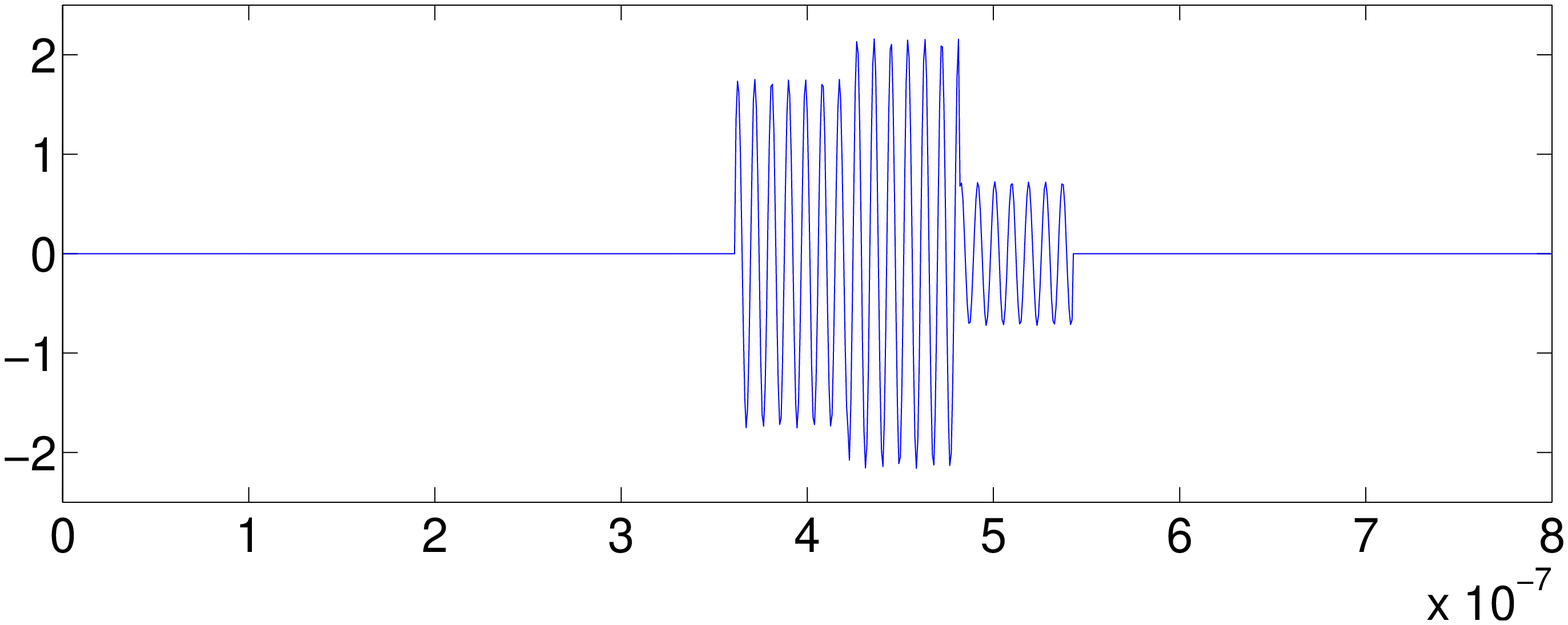}\\
	{\footnotesize(e) reconstructed signal after one iteration}&{\footnotesize(f) reconstructed signal after two iterations}\\
	\includegraphics[scale=.22]{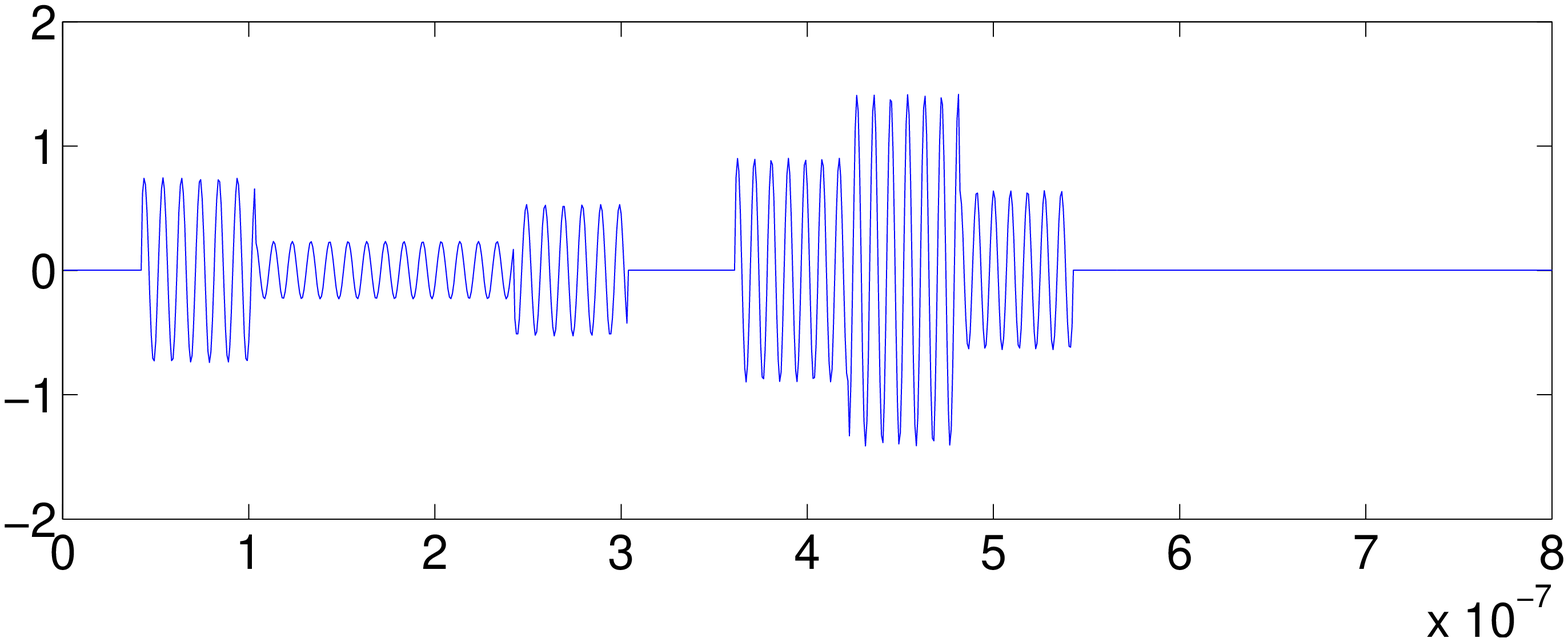}&\includegraphics[scale=.22]{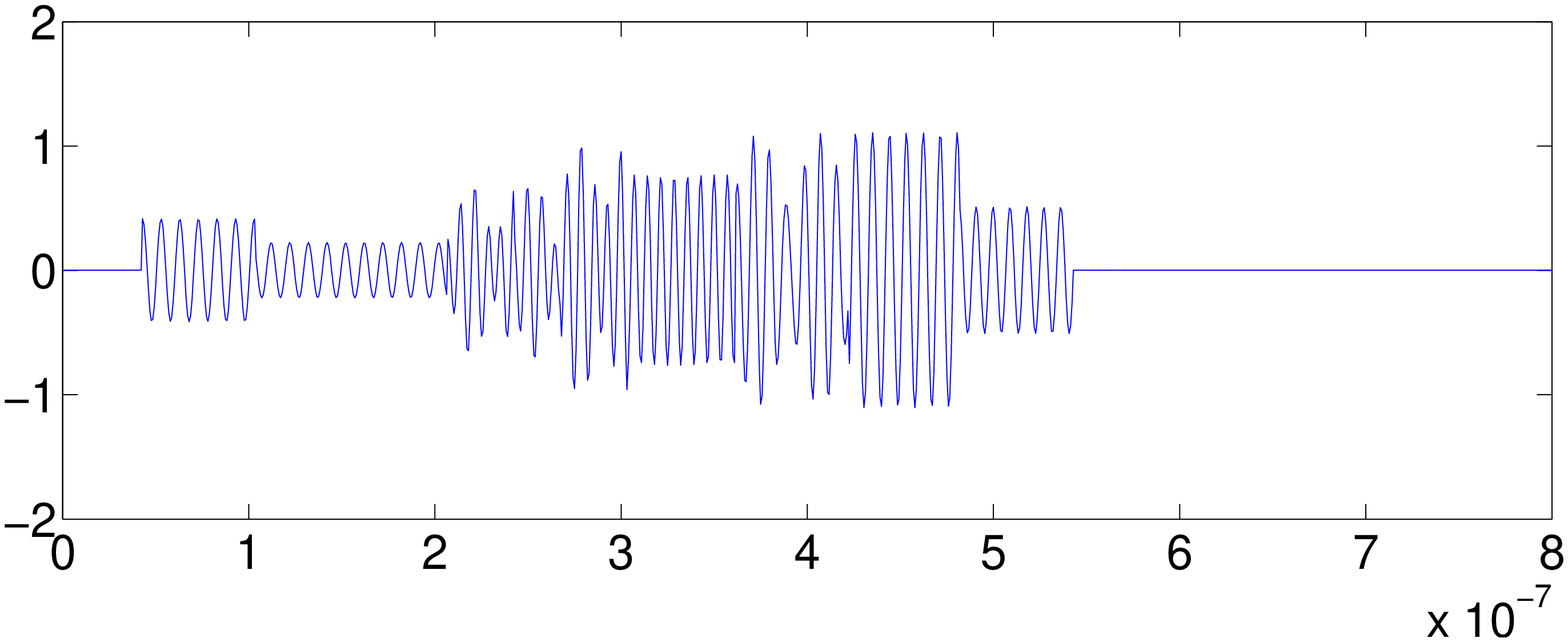}\\
	{\footnotesize(g) reconstructed signal after three  iterations}&{\footnotesize(h) reconstructed signal after four iterations}\\
      \end{tabular}
      \caption{Radar signal and its reconstruction.  All signals, except (a), are on the first sensor. (e)-(h) represent the real part.}
	\label{fig_orig_signal}
	\end{center}
\end{figure}

\subsection{Control of the squared error}
The degree of sparsity in the $\beta$-variables depends on the chosen dictionary. In general, very few prior information is available on the number $\tilde J$ of basis elements used to decompose the radar signal.

A simple but useful tool for model selection consists in studying the evolution of the squared error, that is the evolution of the square norm of the current residual $L(\tilde w^{(t)})=\|Y-\mathbf X\tilde w^{(t)}\|^{2}_{2}$, where  $\tilde w^{(t)}$ denotes the solution at the $t$-th iteration.  Obviously, the squared error $L(\tilde w^{(t)})$ is monotone decreasing over the iterations.  
 
  The OMP-type algorithm always adds the group of variables  that are the most correlated with the current residual, that is the components that are likely to  diminish the squared error the most.  This explains why the algorithm tends to first activating all `true' variables, and afterwards selecting variables that are not in the true model. Hence the problem of selecting the model order  becomes the detection of the point when all `true' variables are activated and the algorithm starts adding pure noise variables.
  
  In fact, activating a `true' variable generally improves the squared error a lot, since the new variable explains a considerable part of the observations. Whereas adding a variable which is not in the true model is less beneficial and  the squared error hardly decreases. 
Consequently,  first the squared error $L(\tilde w^{(t)})$ decreases steeply until the point where all `true' variables are activated. Then the graph of the squared error becomes very flat. This change of the slope yields a kink in the graph of the squared error $L(\tilde w^{(t)})$. Thus, the kink indicates  the number of `true' variables. 

An efficient way to determine the kink consists in computing the convex hull of the graph. Then the kink is given by the endpoint of the longest linear piece  of the convex hull \cite{lavielle05,lavielle05b}.

We may hence proceed as follows. The OMP-type algorithm is run with some fixed $U_\text{max}$ and  stopped  after a large number of iterations. It is useful to store all intermediate solutions $\tilde w^{(t)}$  and the associated square losses.  
Afterwards we determine the final solution as the intermediate solution $\tilde w^{(t)}$  where the squared error presents a kink.

Clearly, the kink is especially pronounced when the signal-to-noise ratio is high.  Conversely, when there is too much noise in the data, then it may be difficult to make out a kink. Besides, when the true model contains a relatively `weak' variable, that is, its contribution to the signal is rather low, then the number of variables  may be underestimated as the `weak' variable can be confounded with noise.

\section{Experimental Results} \label{sec_simul}

\subsection{Setting} 

We illustrate the algorithm by an example. Let us consider a radar signal made of three subsequent sinusoids with different frequencies (100, 140, 110 MHz) and a gap between the two last sinus waves (sinusoids with starting points  0, 150, 300 ns and lengths 150, 100, 120ns), see Figure \ref{fig_orig_signal}(a). The signal is detected by a uniform linear array (ULA) of 5 sensors separated by half a wavelength of the actual narrowband source signals. The sample frequency is $F_s=1.28$ GHz and the total duration of observation is $0.8~\mu$s. Thus, we have $M=1024$ observation points  and $Y$ is a vector of length 5120. Furthermore, we consider three  pathways ($U=3$) with unknown angles of arrivals (9\textdegree, 13\textdegree,  38\textdegree) and delays (60, 95, 120ns). The signal (without noise) received at the first sensor of the antenna is presented in Figure  \ref{fig_orig_signal}(b) and (c), where the superposition of the sinusoids results from the multiple propagation pathways. From these figures it is clear that the multiple paths have a strong impact on the signal and hence cannot be neglected in the analysis. The signal is corrupted by additive noise with a signal-to-noise ratio of 5 dB, see Figure~\ref{fig_orig_signal}(d). The signal-to-noise ratio  is determined on the direct propagation pathway, defined as the 10 $\log$ (energy of the signal on the first path/variance of the noise) as used in \cite{manickam94} or \cite{fuchs99}. 

We choose a Fourier dictionary containing  sine and cosine functions with frequencies 90, 95,\dots, 150 MHz, 
signal lengths 80, 85,\dots, 200 ns and starting points 0, $3\delta$, 6$\delta$,\dots, 255$\delta$ ns, where $\delta=1/F_s\approx0.78$ ns. The size of the dictionary is thus $J= 166\,400$. For  the angles and the delays we use the grids \{1\textdegree, 2\textdegree,\dots, 90\textdegree\} and \{0, $3\delta$, 6$\delta$,\dots, 255$\delta$\}, respectively Hence, $w$ is a vector of dimension $3.83\cdot10^9$.

\subsection{Signal reconstruction by the OMP-type algorithm}
In the simulations we  set $U_\text{max}=2$, that is, at best we recover two of the three pathways and we stop the  algorithm after 9 iterations. Figure \ref{fig_orig_signal}(e)-(h) illustrate the first four iterations of the OMP-type algorithm. In the first iteration the algorithm selects a single sinusoid. Next, a second pathway is selected leading to two overlapping sinusoids.  As now two $\alpha$-variables are activated and $U_\text{max}=2$, the algorithm can only add further $\beta$-variables, that is, further dictionary elements. Thus, in the third iteration  another sinusoid is added which appears directly on the two pathways. Finally, another sinusoid is activated. 

We note that the graph in (h) is quite similar to the one in (c). Nevertheless, the reconstruction is not perfect, as we allowed the algorithm to find only two pathways, while indeed there are three. In this example the second pathway corresponding to the angle of arrival of 13\textdegree~is missing.

\begin{table}[b]{
\begin{center}
\caption{ Empirical mean and standard deviation of the total number of selected variables of $\alpha$- and $\beta$-coefficients  with $U_\text{max}=2$.}\label{tab_squared_error}
\begin{tabular}{rccc}
SNR &{total nb}&{nb of $\alpha$'s }&{nb of $\beta$'s } \\
\hline
5 dB   &5 (0)&  2 (0)&  3 (0) \\
-5  dB  &5 (0)&  2 (0)&  3 (0) \\
-15  dB & 5.1 (0.3) & 2 (0)& 3.1 (0.3) 
\end{tabular}
\end{center}}
\end{table}

\begin{table}[b]{\footnotesize
\begin{center}
\caption{Empirical mean and  standard deviation (in parentheses) of the parameter estimates obtained by the elementary method. }\label{tab_simul_standard}
\begin{tabular}{lll}
\multicolumn{3}{c}{Elementary Method} \\
\hline
&{1st angle (DOA)}&{1st delay (TOA)}\\
\hline
true&{6\textdegree}&{10 ns} \\
\hline
SNR 20 dB&5.99 (0.26) & 10.94 (0) \\
SNR 0 dB &6.81 (8.43)& 22.31 (12.11) \\
\end{tabular}
\end{center}}
\end{table}

\begin{table}[b]{\footnotesize
\begin{center}
\caption{Empirical mean and  standard deviation (in parentheses) of the parameter estimates obtained by the OMP-type algorithm. }\label{tab_simul}
\begin{tabular}{lllll}
\multicolumn{5}{c}{ } \\
\multicolumn{5}{c}{OMP-type algorithm} \\
\hline
&{1st angle}&{1st delay}&{2nd angle}&{2nd delay}\\
\hline
{true}&{6\textdegree}&{10 ns}&{42\textdegree}&{40 ns} \\
\hline
20 dB&6.00 (0) & 10.94 (0) &42.0 (0)&40.6 (0)\\
  0 dB&5.94 (0.24) & 10.5 (1.34) &41.9 (0.27)&40.2 (1.26) \\
\multicolumn{5}{c}{ } \\
&{1st frequ.}&{2nd frequ.}&{3rd frequ.}&{4th frequ.}\\
\hline
{true}&{210 MHz}&{230 MHz}&{210 MHz}&{230 MHz} \\
\hline
20 dB&230.0 (0)& 210.0 (0) &230.0 (0) &210.0 (0) \\
  0 dB&230.0 (0)& 210.0 (0)&230.1 (0.40) &210.0 (0)  \\
\multicolumn{5}{c}{ } \\
&{1st length}&{2nd length}&{3rd length}&{4th length}\\
\hline
{true}&{160 ns}&{160 ns}&{80 ns}&{160 ns} \\
\hline
20 dB&160.0 (0)& 160.6 (1.64) &162.4 (12.1) &170.0 (0) \\
  0 dB&163.1 (4.62)& 165.5 (4.20)&162.4 (15.1) &170.0 (0)  \\
\multicolumn{5}{c}{ } \\
&{1st start}&{2nd start}&{3rd start}&{4th start}\\
\hline
{true}&{0 ns}&{160 ns}&{320 ns}&{400 ns} \\
\hline
20 dB&0 (0)& 158.8 (1.54) &265.3 (12.9) &389.1 (0) \\
  0 dB&0 (0)& 153.1 (5.37) &271.0 (19.7) &389.5 (1.34)  \\
\end{tabular}
\end{center}}
\end{table}

\subsection{Model selection device}
\begin{figure}
	\begin{center}
	\includegraphics[scale=.28]{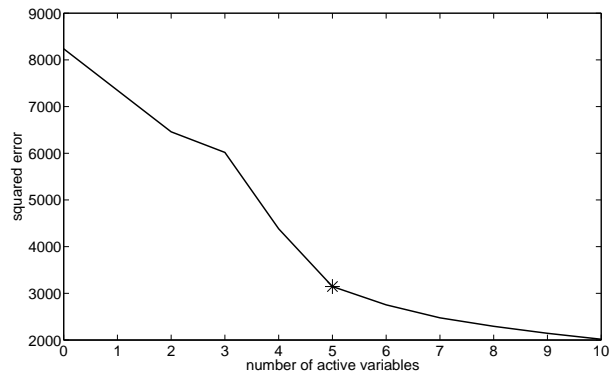}
	\caption{Graph of the squared error of a  dataset with a signal-to-noise ratio of 5 dB. The mark  indicates the kink, that is the selected number of model variables.}
	\label{fig_kink}
	\end{center}
\end{figure}

\begin{figure}
	\begin{center}
	\includegraphics[scale=.26]{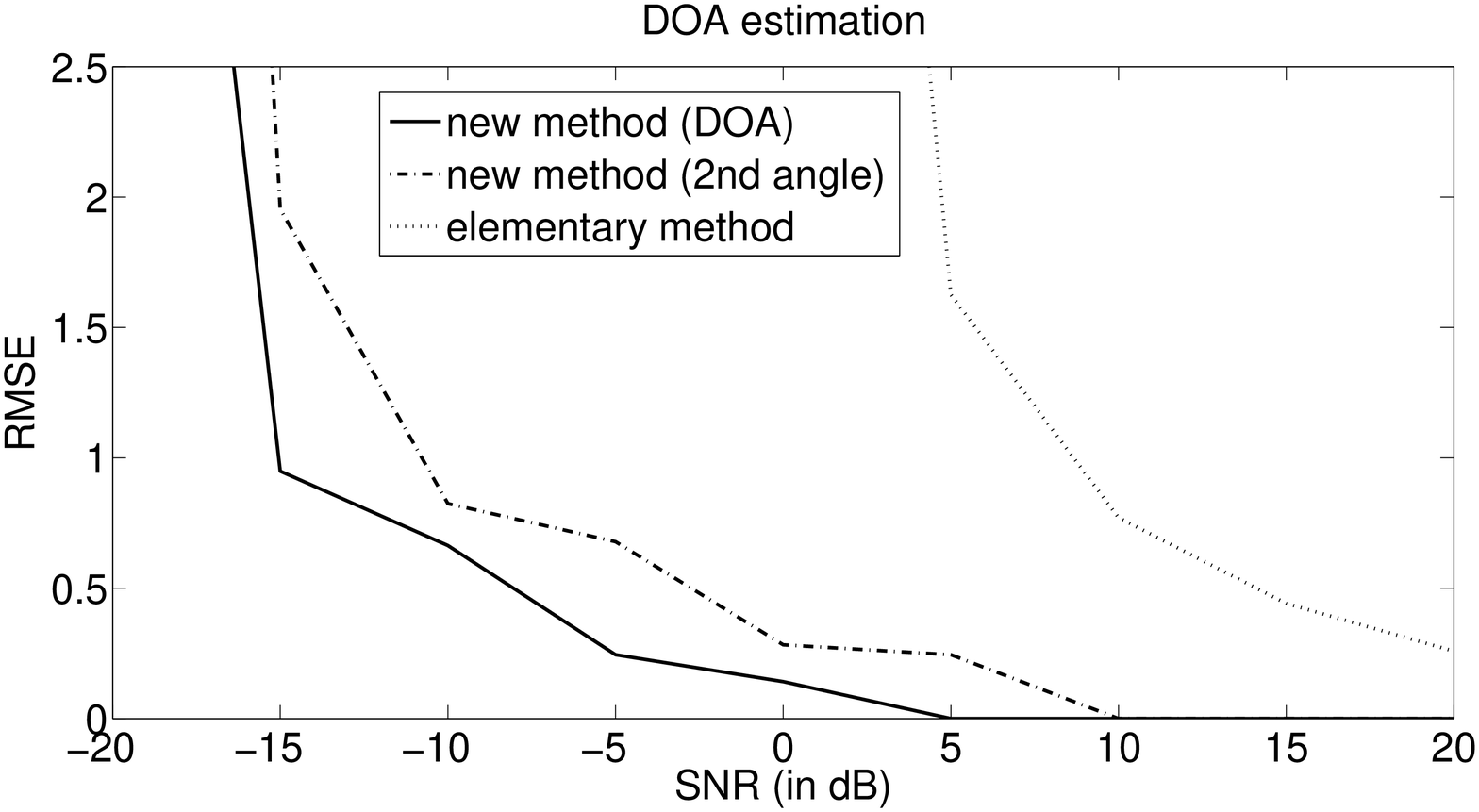}\\
~\\
		\includegraphics[scale=.26]{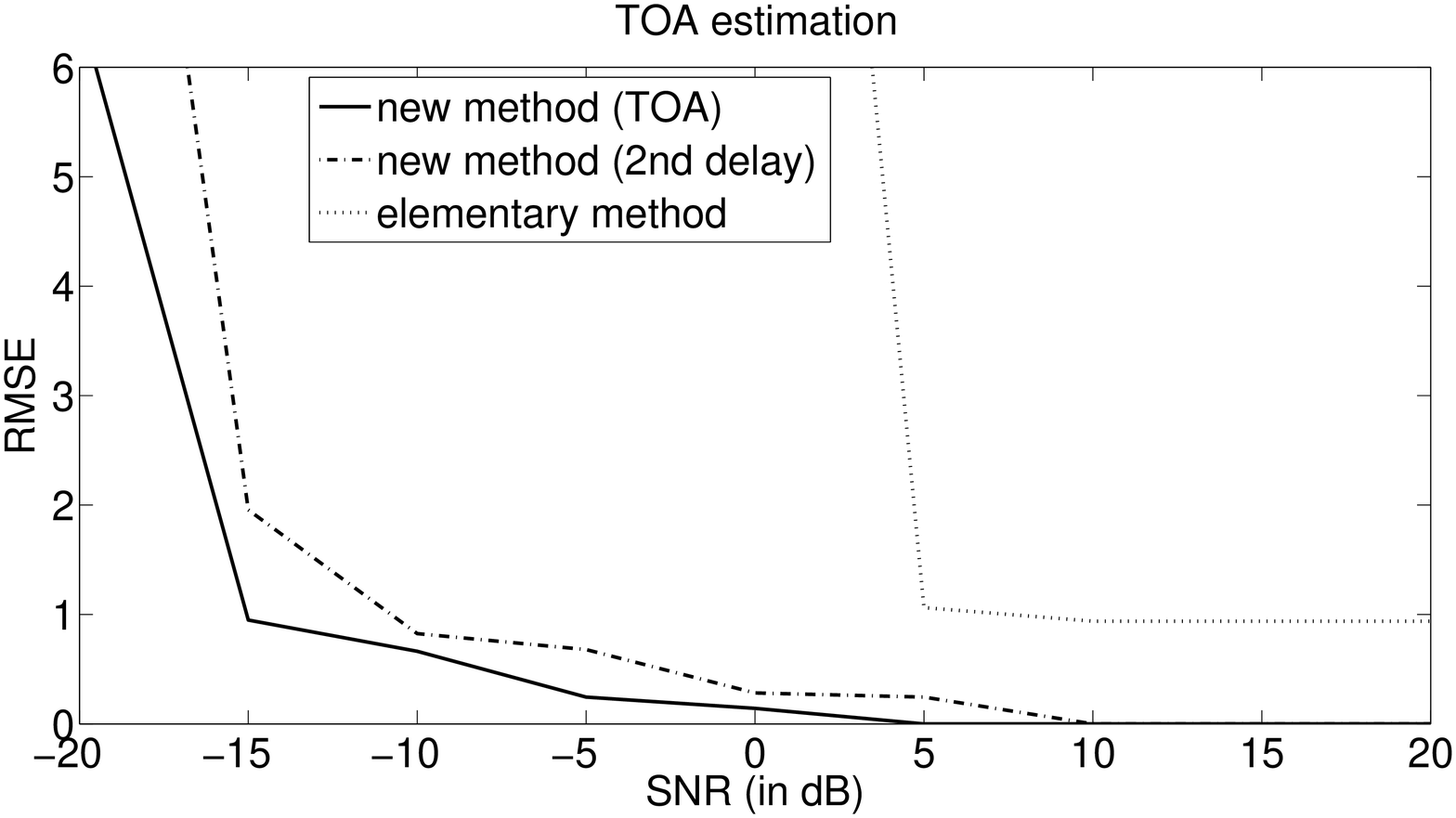}
	\caption{Comparison of RMSE of the DOA and TOA 
	for the elementary method and the new estimator for different signal-to-noise ratios applied on datasets of size 5 x 1280. The RMSE values are based on 50 repetitions per noise level.}
	\label{fig_doa_toa}
	\end{center}
\end{figure}

As mentioned above,  we stopped the algorithm after 9 iterations. Then we used the graph of the squared error to determine the number of variables in the model by the point where the curve exhibits a kink. The kink is determined by computing the convex hull of the graph. 
Figure \ref{fig_kink} illustrates the evolution of the squared loss  for a  sample with a signal-to-noise ratio of 5 dB.  The curve has a slight kink at 5, which is the appropriate number of variables in this case, as the signal is composed of three sinusoids and we set $U_\text{max}=2$. 

 Table  \ref{tab_squared_error} shows  the empirical mean and standard deviation of the total number of selected variables, the number of $\alpha$- and $\beta$-coefficients  for  different signal-to-noise ratios over  50 repetitions. Obviously,   in this example the model selection tool works very well, as the model selection device successfully distinguishes `true' variables from noisy variables, even though  the algorithm is prevented from selecting all `true' variables. Hence, the introduction of an upper bound $U_\text{max}$ for the number of $\alpha$-coefficients does not corrupt the model selection.

\subsection{Comparison with an elementary method} 
A simple, elementary method for estimating the principal DOA in the radar setting consists in first estimating the TOA by thresholding and then  using  the point of observation at the TOA to fit the angle of arrival. More precisely, to detect the TOA by the elementary method, the absolute value of the observed signal at some time point $t$ is compared to some threshold,  here $3\sigma$, where $\sigma^{2}$ is the variance of the noise distribution that is supposed to be known for the elementary method. The first point $t$ where this happens is considered as the TOA. This detection method fails when the signal-to-noise ratio is too low, since then the  signal is too weak compared to the noise such that the observed signal never exceeds the threshold. And when no TOA is detected, the elementary method does not provide an estimate of the DOA.

We compare this elementary method to our OMP-type algorithm in the following setting with an FSK-modulated Barker sequence. We consider the 7-bit Barker sequence [1, 1, -1, -1, 1,-1, -1] with  FSK-modulation. Every bit is represented by a sinusoid of length 80 ns with frequency 210 or 230 MHz. As there are four sign changes in the sequence, the transmitted signal is composed of four  consecutive sinusoids of varying length. Furthermore, we consider two pathways ($U=2$) with unknown angles of arrivals (6\textdegree and 42\textdegree) and delays of 10 and 40 ns.  

Again, the sample frequency is $F_s=1.28$ GHz and we use a ULA of 5 sensors  separated by half a wavelength of the actual narrowband source signals. The observation time  is 1 $\mu$s, such that $Y$ is a vector of length 6400. The dictionary contains the   sine and cosine functions with frequencies 200, 202,\dots, 240 MHz,  lengths 50, 55,\dots, 170 ns and starting points 0, $\delta$, 2$\delta$,\dots, 255$\delta$ ns, where $\delta=1/F_s$. The size of the dictionary is thus $J= 268\,800$.  We use the grids \{1\textdegree, 2\textdegree,\dots, 90\textdegree\} and \{0, $\delta$, 2$\delta$,\dots, 255$\delta$\} for the angles and the delays, respectively Hence, $w$ is a vector of dimension $6.19\cdot10^9$. Data are simulated with different levels of the signal-to-noise ratio. Again we choose $U_\text{max}=2$.

Figure \ref{fig_doa_toa} compares the root mean square error (RMSE) of the DOA and TOA estimates obtained by  both methods based on 50 simulated datasets for each noise level.  Compared to the elementary procedure the OMP-type algorithm achieves a gain of around 20 dB.  This gain is mainly due to the fact that the algorithm takes into account all observations and not only a single point. For low signal-to-noise ratios the elementary method completely fails.

\subsection{Parameter Estimates}

It is important to note that in contrast to the elementary method the OMP-type algorithm provides estimates not only of the DOA and the TOA, but of all the other parameters including the waveform. Indeed, in the setting described above our algorithm identifies several consecutive sinusoids. Table~\ref{tab_simul_standard} presents the empirical mean and standard deviation of the estimates obtained by the elementary method and Table \ref{tab_simul} gives the results for the OMP-type  method, namely for the estimates of the angles of arrival, the associated delays, the frequencies of the four detected sinusoids, their lengths and their starting points. Estimates are especially accurate for the angles of arrival and the frequencies, whereas the estimates of the lengths of the sinusoids and their starting points are less accurate. There is a clear tendency to overestimate the length, especially for the shortest sinusoid.  Indeed, for high noise levels it is not clear if the OMP-type algorithm detects the shortest (and less energetic) sinusoid. The fourth sinusoid generally has not even the correct frequency, so it is possible that the algorithm adds any non relevant  predictor, because the shortest sinusoid is too weak for detection. 

We observe that the reconstructed radar signal $\hat s$ has not exactly the form of  a modulated signal due to the overlapping sinusoids.  It is hence  up to the user to make such an interpretation of the signal. Nevertheless, we emphasize that the algorithm works for any type of modulation (FSK or PSK) without prior knowledge on the modulation and it provides a good idea of the waveform.

\section{Conclusion}\label{sec_concl}
The  algorithm proposed in this paper  is to the best of our knowledge the first method addressing the problem of estimating both the DOA and the unknown waveform of a transmitted radar signal in the presence of multiple propagation pathways. The new method is in the spirit of OMP extended to more sophisticated sparsity structures.
The scaling of the OMP-type algorithm has been considered, which is necessary for unbiased estimation. A simulation study  illustrates the good performance of the new method and exhibits a considerable improvement with respect to some elementary method.

We would like to emphasize the general character of our method, which is of much interest for the intercepted radar signal setting. First, the method proposed here is not restricted to  a specific  array geometry.  Second and more importantly, the approach can be adapted to a large variety of signals such as sinusoids, chirps, composed signals etc. by using a dictionary where a sparse representation of the signal is possible. 

\bibliographystyle{IEEEtran}

\bibliography{biblio}

\end{document}